\documentclass[letterpaper,12pt,twoside]{article}

\usepackage[english]{babel}
\usepackage{times}
\usepackage{epsfig}
\usepackage{amsfonts}
\usepackage{amsthm}
\usepackage{float}
\usepackage{amsmath}
\usepackage{amssymb}
\usepackage{latexsym}
\usepackage{graphicx, color}
\usepackage{caption}
\usepackage{subcaption}
\usepackage{prodint}
\usepackage[round]{natbib}
\usepackage{bm}
\usepackage{textcomp}
 \usepackage{paralist}
\usepackage{tikz}
\usepackage{mathabx}
\usepackage{setspace}
\usetikzlibrary{arrows,automata,positioning}

\newcommand{\rnc}{\renewcommand}
\newcommand{\nc}{\newcommand}
\newcommand{\mrm}{\mathrm}
\renewcommand{\b}{\textbf}

\nc{\mb}{\mathbb}
\nc{\mac}{\mathcal}
\nc{\E}{\mb{E}}
\nc{\N}{\mb{N}}
\nc{\R}{\mb{R}}
\nc{\Q}{\mb{Q}}
\rnc{\P}{\mrm P}
\rnc{\d}{\mrm d}
\nc{\C}{\mac{C}}
\nc{\D}{\mac{D}}
\nc{\B}{\mac{B}}
\nc{\wh}{\widehat}
\nc{\oPo}{\stackrel{p}{\longrightarrow}}
\nc{\oWo}{\stackrel{w}{\longrightarrow}}
\nc{\oDo}{\stackrel{d}{\longrightarrow}}

\relax 
\allowdisplaybreaks[3]

\relax

\newcommand{\bay}{\begin{array}}
\newcommand{\eay}{\end{array}}

\newcommand{\bqa}{\begin{eqnarray*}}
\newcommand{\eqa}{\end{eqnarray*}}

\newcommand{\bee}{\begin{eqnarray*}}
\newcommand{\eee}{\end{eqnarray*}}

\newcommand{\bea}{\begin{eqnarray*}}
\newcommand{\eea}{\end{eqnarray*}}

\newcommand{\bqan}{\begin{eqnarray}}
\newcommand{\eqan}{\end{eqnarray}}

\newcommand{\be}{\begin{eqnarray}}
\newcommand{\ee}{\end{eqnarray}}

\newcommand{\bit}{\begin{itemize}}
\newcommand{\eit}{\end{itemize}}

\newcommand{\ben}{\begin{enumerate}}
\newcommand{\een}{\end{enumerate}}

\newcommand{\beq}{\begin{equation}}
\newcommand{\eeq}{\end{equation}}

\newcommand{\bdes}{\begin{description}}
\newcommand{\edes}{\end{description}}

\newcommand{\btb}{\begin{tabular}}
\newcommand{\etb}{\end{tabular}}

\newcommand{\bcen}{\begin{center}}
\newcommand{\ecen}{\end{center}}

\newcommand{\bmp}{\begin{minipage}}
\newcommand{\emp}{\end{minipage}}

\newcommand{\vX}{\boldsymbol{X}}

\newcommand{\vZ}{\boldsymbol{Z}}

\newtheorem{definition}{{\sc Definition}\sc}[section]
\newcommand{\bdefi}{\begin{definition}}
\newcommand{\edefi}{\end{definition}}

\newtheorem{appropr}[definition]{{\sc Approximation Procedure}\sc}
\newcommand{\bappr}{\begin{appropr}}
\newcommand{\eappr}{\end{appropr}}

\newtheorem{bedi}[definition]{{\sc Condition}\sc}
\newcommand{\bbd}{\begin{bedi}}
\newcommand{\ebd}{\end{bedi}}

\newtheorem{bedin}[definition]{{\sc Conditions}\sc}
\newcommand{\bbdn}{\begin{bedin}}
\newcommand{\ebdn}{\end{bedin}}

\newtheorem{corollary}[definition]{{\sc Corollary}\sc}
\newcommand{\bco}{\begin{corollary}}
\newcommand{\eco}{\end{corollary}}

\newtheorem{proposition}[definition]{{\sc Proposition}\sc}
\newcommand{\bpro}{\begin{proposition}}
\newcommand{\epro}{\end{proposition}}

\newtheorem{satz}[definition]{{\sc Theorem}\sc}
\newcommand{\bsa}{\begin{satz}}
\newcommand{\esa}{\end{satz}}

\newtheorem{assumption}[definition]{{\sc Assumption}\sc}
\newcommand{\bas}{\begin{assumption}}
\newcommand{\eas}{\end{assumption}}

\newtheorem{assumptions}[definition]{{\sc Assumptions}\sc}
\newcommand{\bass}{\begin{assumptions}}
\newcommand{\eass}{\end{assumptions}}

\newtheorem{abb}{{\sc Figure}\sc}
\newcommand{\babb}{\begin{abb}}
\newcommand{\eabb}{\end{abb}}

\newenvironment{remark}{\begin{rmk}\sl}{\end{rmk}}
\newtheorem{rmk}{{\sc Remark}\sc}[section]
\newcommand{\brem}{\begin{remark}}
\newcommand{\erem}{\end{remark}}

\newenvironment{remarks}{\begin{rmks}\sl}{\end{rmks}}
\newtheorem{rmks}{{\sc Remarks}\sc}[section]
\newcommand{\brems}{\begin{remarks}}
\newcommand{\erems}{\end{remarks}}

\newenvironment{example}{\begin{exmp}\rm}{\end{exmp}}
\newtheorem{exmp}{{\sc Example}\sc}[section]
\newcommand{\bbsp}{\begin{example}}
\newcommand{\ebsp}{\end{example}}
\newcommand{\bexa}{\begin{example}}
\newcommand{\eexa}{\end{example}}

\newtheorem{model}{{\sc Model}\sc}[section]
\newcommand{\bmdl}{\begin{model}}
\newcommand{\emdl}{\end{model}}

\newtheorem{scheme}{{\sc Scheme}\sc}[section]
\newcommand{\bscm}{\begin{scheme}}
\newcommand{\escm}{\end{scheme}}

\newenvironment{tabelle}{\begin{tabl}\sl}{\end{tabl}}
\newtheorem{tabl}{{\sc Table}\sc}
\newcommand{\btab}{\begin{tabelle}}
\newcommand{\etab}{\end{tabelle}}

\newenvironment{exercise}{\begin{exc}\sl}{\end{exc}}
\newtheorem{exc}{Exercise}[section]
\newcommand{\bexe}{\begin{exercise}}
\newcommand{\eexe}{\end{exercise}}

\numberwithin{equation}{section}

\newtheorem{thm}[definition]{Theorem}
\newtheorem{lemma}[definition]{Lemma}
\newtheorem{rem}[definition]{Remark}

\begin{document}

\linespread{1.3}

    
\title{\Large\bf Bootstrap- and Permutation-Based Inference \\ for the Mann-Whitney Effect \\ for Right-Censored and Tied Data
}
    
\author{Dennis Dobler$^{*}$ and Markus Pauly \\[1ex] 
}

\date{}

\maketitle

\vspace{-1cm}

{\centering
{\small Ulm University, \ \
Institute of Statistics \\
Helmholtzstr. 20, \ \
89081 Ulm, \ \
Germany \\
email: dennis.dobler@uni-ulm.de, \ \
email: markus.pauly@uni-ulm.de \\
}}

\vspace{0.2cm}

{
\centering

\date{\today}

}

\begin{center}
{\bf Summary} \vspace{-.4cm}
\end{center}
\noindent The Mann-Whitney effect is an intuitive measure for discriminating two survival distributions. 
 Here we analyze various inference techniques for this parameter in a two-sample survival setting with independent right-censoring, 
 where the survival times are even allowed to be discretely distributed. This allows for ties in the data and requires 
 the introduction of normalized versions of Kaplan-Meier estimators from which adequate point estimates are deduced.
 {\color{black}Asymptotically exact inference procedures based on standard normal, 
 bootstrap- and permutation-quantiles are developed and compared in simulations.} 
 Here, the asymptotically robust and {\color{black} -- under exchangeable data -- }
 even finitely exact permutation procedure turned out to be the best.
 Finally, all procedures are illustrated using a real data set.

\noindent{\bf Keywords:} 
Counting process; 
Efron's bootstrap;
Heteroscedasticity;
Kaplan-Meier estimator; 
Permutation technique.

\vfill
\vfill

\vspace{0cm}

\noindent${}^{*}$ {Corresponding author}

\pagenumbering{gobble}
\thispagestyle{empty}

\newpage

\pagenumbering{arabic}

\allowdisplaybreaks

\section{Introduction}

When comparing the survival times from two independent groups ($j=1,2$) the {\it Mann-Whitney effect}
is an intuitive measure; see e.g. \cite{koziol2009concordance}.
In a classical survival setting with continuous life time distributions and random censoring it is given by the probability
$ P(T_{1} > T_{2}) $
that a random subject from group $j=1$ (with survival time $T_1$) survives longer than a randomly chosen person from group $j=2$ (with survival time $T_2$). 
In case of uncensored data this effect reduces to the well-known Wilcoxon functional underlying the nonparametric Mann-Whitney test. Depending on the field of application 
it is also known as {\it nonparametric treatment effect} (e.g. \citealt{brunner2000nonparametric}), {\it stress-strength characteristic} (e.g. \citealt{kotz2003stress}) or {\it probabilistic index} (e.g. \citealt{thas12}). 
Moreover, in case of diagnostic tests it has a direct interpretation as the area under the corresponding ROC curve, see e.g. 
\cite{lange2012sensitivity}, \cite{pauly2016} and \cite{zapf2015wild} as well as \cite{zhou2009statistical} for more details on diagnostic accuracy measures. 
The Mann-Whitney effect is often estimated by the $c$-index for concordance (e.g. \citealt{koziol2009concordance}).
As pointed out by \cite{Acion2006}, 
this is ``{\it a simple, clinically relevant, and robust index}'' and thus ``{\it an ideal effect measure}'', see also \cite{Kieser2013}. 
The same still holds true in case of survival outcomes that may be subject to independent random censoring, see e.g. the glorification of the $c$-index in \cite{hess2010comparing} or \cite{dunkler2010gene}. 
An {\tt R}-package for a Wilcoxon-Mann-Whitney test was propagated in \cite{de2014unifiedwmwqpcr}.

In the present paper we face the practically relevant situation where tied data are often inevitable. 
Thus, in order to take ties appropriately into account, we use a more general definition of the Mann-Whitney effect:
\begin{equation}\label{eq: C-Index}
 p=P(T_{1} > T_{2}) + \frac12 P(T_{1} = T_{2}),
\end{equation}
also known as {\it ordinal effect size measure} in case of completely observed data (\citealt{ryu2008modeling, konietschke2012rank}). 
Recently, a related effect measure, the so called {\it win ratio} (for prioritized outcomes), has been investigated considerably by several authors 
(\citealt{pocock2012win}, \citealt{rauch2014opportunities}, \citealt{luo2014alternative}, \citealt{abdalla2016win, bebu16} as well as \citealt{wang2016win}). 
It is the odds of the Mann-Whitney effect $p$ in  \eqref{eq: C-Index}, i.e.
\begin{equation}\label{eq: win ratio}
  w = \frac{p}{1-p},
\end{equation}
which is also referred to as the {\it odds of concordance}, see \cite{martinussen2013estimation} for a treatment in the context of a semiparametric regression model. 
In our situation $p$ and $w$ describe the probability that a patient of group 1 survives longer than a patient of group 2.
That is, $p > 1/2$, or equivalently $w > 1$, imply a protective survival effect for group 1.
Note that until now, ties have been excluded for estimating these quantities which particularly led to the recent assessment of 
\cite{wang2016win} that ``{\it we caution that the win ratio method should be used only when the amount of tied data is negligible}''.

In this paper, we propose and rigorously study different statistical inference procedures for both parameters $p$ and $w$ in a classical survival 
model with independent random censoring, even allowing for ties in the data. 
While several authors (e.g. \citealt{nandi1994note}, \citealt{cramer1997umvue}, \citealt{kotz2003stress}, and references therein) have considered inference for $p$
under specific distributional assumptions,
we here focus on a completely nonparametric approach, not even assuming continuity of the data. 
Apart from confidence intervals for $p$ and $w$ this also includes one- and two-sided test procedures for the null hypothesis of no group effect (tendency)
\begin{equation}\label{eq: null}
H_0^p: \Big\{ p =  \frac12 \Big\}  = \Big\{ w =  1 \Big\}.
\end{equation}
In the uncensored case this is also called the {\it nonparametric Behrens-Fisher problem}, see e.g. \cite{brunner2000nonparametric} and 
\cite{neubert07}. To this end, the unknown parameters $p$ and $w$ are estimated by means of normalized versions of Kaplan-Meier estimates. 
These are indeed their corresponding nonparametric maximum likelihood estimates, see \cite{efron67} as well as \cite{koziol2009concordance} for the case of continuous observations. 
Based on their asymptotic properties we derive asymptotically valid tests and confidence intervals. 
These may be regarded as extensions of the Brunner-Munzel (2000) test to the censored data case. 
Since, for small sample sizes, the corresponding tests may lead to invalid $\alpha$-level control 
(e.g. \citealt{medina2010inappropriate} or \citealt{pauly2016} without censoring) we especially discuss and analyze two different 
resampling approaches (bootstrapping and permuting) to obtain better small sample performances. 

{
The resulting tests are innovative in several directions compared to other existing procedures for the two-sample survival set-up:
\begin{enumerate}
 \item We focus on the null hypothesis $H_0^p$ in \eqref{eq: null} of actual interest.
 Before, only the more special null hypothesis $H_0^S : \{ S_1 = S_2 \}$ of equal survival distributions between the two groups has been investigated, 
  see e.g. \cite{efron67}, \cite{akritas97} and \cite{akritas2011nonparametric}. 
  Corresponding one-sided testing problems (for null hypotheses formulated in terms of distribution functions) 
  based on the related {\it stochastic precedence} have been treated in \cite{arcones02} and \cite{davidov2012ordinal}. 
  Instead, our procedures are not only able to assess the similarity of two survival distributions but they also quantify the degree of deviation 
  by confidence intervals for meaningful parameters. 
 \item The more complex null $H_0^p$ has so far only been studied in the uncensored case, 
  see e.g. \cite{janssen1999}, \cite{brunner2000nonparametric}, \cite{de2013extension}, \cite{chung2016asymptotically}, \cite{pauly2016} and the references cited therein. 
  The present combination of allowing for survival analytic complications and focussing on the effect size $p$ 
  is achieved by utilizing empirical process theory applied to appropriate functionals.
 \item We do not rely on the (elsewhere omnipresent) assumption of existing hazard rates.
  Instead, we here adjust for tied data by considering normalized versions of the survival function and the Kaplan-Meier estimator (leading to mid-ranks in the uncensored case). 
  This more realistic assumption of ties in the observations accounts for a phenomenon
  which is oftentimes an intrinsic problem by study design.
  Therefore, methodology for continuous data (even for only testing $H_0^S$) should not be applied. 
  Notable exceptions for the combination of survival methods and discontinuous data are provided in \cite{akritas97} and \cite{brendel2014weighted} where the hypothesis $H_0^S$ is tested.
 \item Finally, small sample properties of inference procedures relying on the asymptotic theory are greatly improved by applications of resampling techniques.
 These utilized resampling techniques are shown to yield consistent results even in the case of ties.
  Thereof, the permutation procedure succeeds in being even finitely exact in the case of exchangeable survival data in both sample groups; see e.g. \cite{lehmann06}, \cite{good10}, \cite{pesarin10,pesarin12}, and \cite{bonnini14} for the classical theory on permutation tests.
\end{enumerate}
}
In this perspective, the present paper not only states the first natural extension of point estimates for $p$ to tied survival data,
but especially introduces the first inference procedures for $H_0^p$ (tests and confidence intervals) with rigorous consistency proofs. 
The latter have formerly not even been known in the continuous survival case.


The article is organized as follows. Section~\ref{sec:model} introduces all required notation and estimators,
whose combination with the variance estimator in Section~\ref{sec:stud} yields (non-resampling) inference procedures.
Theoretical results concerning the resampling techniques will be presented in Sections~\ref{sec:ebs} (the pooled bootstrap) and~\ref{sec:perm} (random permutation).
A simulation study in Section~\ref{sec:simus} reports the 
improvement of the level $\alpha$ control by the proposed permutation and bootstrap techniques.
A final application of the developed methodology to a tongue cancer data set (\citealt{klein03}) is presented in Section~\ref{sec:example}.
This article's results are discussed in Section~\ref{sec:dis} and theoretically proven in Appendix~\ref{sec:proofs}.

\section{Notation, model, and estimators}
\label{sec:model}

For a more formal introduction of the concordance index $p$ and the win ratio $w$ we employ some standard notation from survival analysis. 
Thus, we consider two independent groups $(j=1,2)$ of independent survival times 
$
\tilde T_{j1}, \dots, \tilde T_{jn_j}: (\Omega,\mathcal{A},P) \rightarrow (0,\infty),\; j=1,2,
$
with total sample size $n = n_1 + n_2, \ n_1, n_2 \in \N$. 
Their distributions may have discrete components which reflects the situation in most clinical studies (i.e. survival times rounded to days or weeks).
Since most studies pre-specify a point of time $K > 0$ after which no further observation is intended,
we may also truncate the above survival times to
$$ T_{ji} = \tilde T_{ji} \wedge K = \min(\tilde T_{ji}, K); \quad i=1,\dots, n_j, \ j=1,2.$$
Denote their survival functions as $S_j(t) = 1 - F_j(t) \equiv P(T_{j1} > t)$, $j=1,2$.
Thus, both sample groups may have different, even heteroscedastic distributions.
Their cumulative hazard functions are given by 
$\Lambda_j(t) = - \int_0^t \frac{\d S_j}{S_{j-}}$ where the index minus (here in $S_{j-}$) always indicates the left-continuous version of a right-continuous function.

The survival times are randomly right-censored by independent, positive variables
$
C_{j1}, $ $\dots,$ $ C_{jn_j} $  
with possibly discontinuous censoring survival functions $G_j, j=1,2$. 
Observation is thus restricted to
\begin{align*}
 \b X_j = \{(X_{j1}, \delta_{j1}), \dots, (X_{jn_j}, \delta_{jn_j})\}, \quad j=1,2,
\end{align*}
where $X_{ji} = \min(T_{ji}, C_{ji}), \delta_{ji} = \b 1\{ X_{ji} = T_{ji} \}, 1 \leq i \leq n_j$.
Note that the choice of $K$ shall imply a positive probability of each $\{\tilde T_{ji} > K\}$.
This constant $K$ could, for example, be the end-of-study time, i.e. the largest censoring time.
For later use we also introduce the usual counting process notation
\begin{align*}
 \!\! N_{j;i}(u) = \ &  \b 1\{ \text{``ind. $i$ of group $j$ has an observed event during $[0,u]$''} \} &&  \!\!\!\!\! =  \b 1 \{ X_{ji} \leq u, \delta_{ij} = 1 \},\\
 \!\! Y_{j;i}(u) = \ &  \b 1\{ \text{``ind. $i$ of group $j$ is under observation at time $u-$''} \}  && \!\!\!\!\! =  \b 1 \{ X_{ji} \geq u \}.
\end{align*}
Summing up these quantities within each group results in $Y_j(u) = \sum_{i = 1}^{n_j} Y_{j;i}(u)$, the number of group-$j$ subjects under study shortly before $u$,
and $N_{j}(u) = \sum_{i=1}^{n_j} N_{j;i}(u)$, the number of observed events in group $j$ until time $u$.
Denote by $f^\pm = \frac12 (f + f_- ) $ the normalized version of a right-continuous function $f$.
With this notation the Mann-Whitney effect \eqref{eq: C-Index} and the win ratio \eqref{eq: win ratio} are given as
\begin{equation}\label{eq: C index detailed}
  p = P(T_{11} > T_{21}) + \frac12 P(T_{11} = T_{21})  = - \int S_1^\pm \d S_2 = 1 - \int F_1^\pm \d F_2
\end{equation}
and $w= {p}/{(1-p)}$, respectively.
If not specified, integration is over $[0,K]$.
In this set-up we test the null hypothesis $H_0^p: \{ p =  1/2 \}  = \{ w =  1 \}$ that the survival times from both groups are 
{\it tendentiously equal} against one- or two-sided alternatives. 
We note that the usually considered null hypothesis $H_0^S:\{S_1 = S_2\}$ of equal survival distributions is more restrictive and implies $H_0^p$. 
Similarly, a stochastic order or precedence (\citealt{davidov2012ordinal}) such as $F_1 \lneq F_2$ implies $p\geq 1/2$.

In order to test $H_0^S$ for continuous survival times, \cite{efron67} has introduced a natural estimator for $p$, see also \cite{koziol2009concordance}, by replacing the unknown survival functions 
$S_j$ in \eqref{eq: C index detailed} with the corresponding Kaplan-Meier estimators $\wh S_j, j=1,2$.  
Thus, utilizing the normalized versions of the non-parametric maximum likelihood estimators yields %
\begin{equation}\label{eq: C index estimate}
 \wh p = \wh P(T_{11} > T_{21}) + \frac12 \wh P(T_{11} = T_{21})  = - \int \wh S_1^\pm \d \wh S_2
\end{equation}
and $\wh w= \wh p / (1-\wh p)$ as reasonable estimators of $p$ and $w$, respectively. 
Similar estimators for $p$ have been proposed by \cite{akritas97} 
and \cite{brunner2000nonparametric}. 
The latter quantity $\wh w$ has been introduced by \cite{pocock2012win} for uncensored observations (without ties) 
with the nice interpretation as total number of winners divided by the total number of losers in group $1$ 
(where $T_{1i}$ wins against $T_{2\ell}$ if $T_{1i} > T_{2\ell}$). 

In order to obtain tests for $H_0^p$ as well as one- or two-sided confidence intervals for the Mann-Whitney effect $p$ and the win ratio $w$, 
we study the limit behaviour of $\wh p$ under the general asymptotic framework
\begin{align}
 \label{cond:liminfsup}
  0 < \liminf \frac{n_1}{n} \leq \limsup \frac{n_1}{n} < 1
\end{align}
as $n_1 \wedge n_2 \rightarrow \infty$.
For a better illustration in intermediate steps we sometimes also assume that a unique limit exists, i.e., as $n_1 \wedge n_2 \rightarrow \infty$,
\begin{align}
\label{cond:liminfsup2}
 n_1/n \rightarrow \kappa \in (0,1).
\end{align}
Our main theorems, however, will only rely on the weak assumption of display~\eqref{cond:liminfsup}.

We  denote by ``$\oDo$'' and ``$\oPo$'' weak convergence and convergence in outer probability 
as $n \rightarrow \infty$, respectively, both in the sense of \cite{vaart96}. 
The following central limit theorem for $\wh p \ $ is the normalized counterpart of the asymptotics due to \cite{efron67} 
and is proven by means of the functional $\delta$-method in combination with the weak convergence theorem for the Kaplan-Meier estimator 
as stated in \cite{vaart96}, Section~3.9.
\begin{thm}
\label{thm:weak} 
\noindent(i) 
Suppose \eqref{cond:liminfsup2} holds, then the Mann-Whitney statistic $V_n=V_n(p) = \sqrt{{n_1 n_2}/{n}} (\wh p - p)$ is asymptotically normal distributed, i.e.  
 $V_n\oDo Z \sim N(0, \sigma^2)$ as $n \rightarrow \infty$.
 The limit variance is $\sigma^2 = (1-\kappa) \sigma^2_{12} + \kappa \sigma^2_{21}$, where for $j,k\in\{(1,2), (2,1)\}$:
 \begin{equation}\label{eq: limit variance}
  \sigma^2_{jk} = \int \int \Gamma_j^{\pm \pm}(u,v)  \d S_k(u) \d S_k(v) 
  \ \text{ and } \
  \Gamma_j(u,v) = S_j(u) S_j(v) \int_0^{u \wedge v} \frac{\d \Lambda_j}{S_{j-} G_{j-}}.
 \end{equation}
 Moreover, $\Gamma_j^{\pm \pm}$ denotes the normalized covariance function given by
 $$ \Gamma_j^{\pm \pm}(u,v) = \frac14 [\Gamma_j(u,v) + \Gamma_j(u-,v) + \Gamma_j(u,v-) + \Gamma_j(u-,v-)]. $$
 \noindent(ii) Under the conditions and notation of (i), the win ratio statistic $ \sqrt{{n_1 n_2}/{n}} (\wh w - w)$ 
 asymptotically follows a normal-$N(0, \sigma^2/(1-p)^4)$-distribution.
\end{thm}

\begin{rem}
\label{rem:weak}
 (a) \cite{efron67} originally used a version of the Kaplan-Meier estimator that always considered the last observation as being uncensored. 
Moreover, he considered the more special null hypothesis $H_0^S : \{ S_1 = S_2 \}$ in order to simplify the variance representation under the null.\\
(b) Without a restriction at some point of time $K$, a consistent estimator for the Mann-Whitney effect requires 
consistent estimators for the survival functions on their whole support.
This involves the condition $-\int (\d S_j) / G_{j-} < \infty, j=1,2,$ which is obviously only possible if the support of $S_j$ is contained in the support of $G_j$;
see e.g. \cite{gill83}, \cite{ying89} and \cite{akritas97}.
However, this assumption is often not met in practice, e.g. if $G_j(u) = 0 < S_j(u)$
for some point of time $u > 0$.
\end{rem}

Since the variances $\sigma^2_{jk}$ are unknown under the null $H_0^p$, 
their estimation from the data is mandatory in order to obtain 
asymptotically consistent tests and confidence intervals for $p$ and $w$.

\section{Variance estimation and studentized test statistics}
\label{sec:stud}

Asymptotically pivotal test statistics result from studentized versions of $\wh p$ and $\wh w$. 
A natural estimator for the limit variance $\sigma^2$ of $V_n= \sqrt{{n_1 n_2}/{n}} ( \wh p - p)$ 
is $\wh \sigma^2 = \frac{n_1 n_2}{n} ( \wh \sigma^2_{12} + \wh \sigma^2_{21} )$, where
\begin{equation}\label{eq: var estimate}
 n_j \wh \sigma^2_{jk} = \int \int n_j \wh \Gamma_j^{\pm\pm}(u,v) \d \wh S_k(u) \d \wh S_k(v) \ \ \text{and} \ \
 \wh \Gamma_j(u,v) = \wh S_j(u) \wh S_j(v) \int_0^{u \wedge v} \frac{\d N_j}{(1 - \frac{\Delta N_j}{Y_j})Y^2_j} 
\end{equation}
for $1 \leq j \neq k \leq 2$.
Here, the function $\Delta f(u) = f(u) - f(u-)$ contains all jump heights of a right-continuous function $f$.
\begin{lemma}
 \label{lem:cons_var_est}
Under \eqref{cond:liminfsup2} the estimator $\wh \sigma^2$ is consistent for $\sigma^2$ defined in Theorem~\ref{thm:weak}, i.e. $\wh \sigma^2 \oPo \sigma^2$.
\end{lemma}
This result directly leads to the studentized statistics
\begin{equation}
 T_n(p) =  \sqrt{\frac{n_1 n_2}{n}} \frac{\wh p - p}{\wh \sigma} \quad \text{ and } \quad
 W_n(w) =  \sqrt{\frac{n_1 n_2}{n}} \frac{(1-\wh p)^2}{\wh \sigma} (\wh w - w) =  \sqrt{\frac{n_1 n_2}{n}} \frac{\wh w - w}{\wh \sigma (1+\wh w)^2} 
\end{equation}
which are both asymptotically standard normal as $n_1 \wedge n_2 \to \infty$ by Slutzky's Theorem and the $\delta$-method only assuming \eqref{cond:liminfsup}.
Indeed, under \eqref{cond:liminfsup}, Theorem~\ref{thm:weak} and Lemma~\ref{lem:cons_var_est} might be applied along each convergent subsequence of $n_1/n$.
Since all resulting limit distributions of $T_n(p)$ and $W_n(p)$ are pivotal, i.e. independent of $\kappa$, this weak convergence must hold for the original sequence as well.
Thus, two-sided confidence intervals for $p$ and $w$ of 
asymptotic level  $(1-\alpha) \in (0,1)$ are
\begin{align}
 \label{eq:consistent_cis}
  {I}_n = \left[\wh p \ \mp \ \frac{z_{\alpha/2}\wh \sigma \sqrt{n}}{\sqrt{n_1n_2}}  \right] \text{ (for $p$)} \quad \text{ and } \quad 
  \left[\wh w \ \mp \ \frac{z_{\alpha/2}\wh \sigma}{(1- \wh p)^2} \sqrt{\frac{n}{n_1n_2}}  \right] \text{ (for $w$)},
\end{align}
respectively, where $z_\alpha$ denotes the $(1-\alpha)$-quantile of the standard normal distribution. 
Moreover, 
\begin{equation}
\label{eq:consistent_test}
\varphi_n = \mathbf{1}\{T_n({1}/{2}) > z_\alpha\}  \quad \text{ and } \quad \psi_n = \mathbf{1}\{ W_n(1) > z_\alpha\}
\end{equation}
are consistent asymptotic level $\alpha$ tests for $H_0^p:\{p=\frac{1}{2}\} = \{ w = 1 \}$ against the one-sided alternative hypothesis $H_1^p:\{p>\frac{1}{2}\} = \{ w > 1 \}$, i.e. 
$E(\varphi_n) \rightarrow \alpha \mathbf{1}\{p=1/2\}  + \mathbf{1}\{p>1/2\}$ and $E(\psi_n) \rightarrow \alpha \mathbf{1}\{w=1\}  + \mathbf{1}\{w>1\}$ as $n\to\infty$. One-sided confidence intervals and two-sided tests can be obtained by inverting the above procedures.
For larger sample sizes ($n_j > 30$ depending on the magnitude of censoring),
the above inference methods \eqref{eq:consistent_cis} and \eqref{eq:consistent_test} are fairly accurate; 
see the simulation results in Section~\ref{sec:simus}. 
For smaller sample sizes, however, these procedures tend to have inflated type-$I$ error probabilities. 
Therefore, we propose different resampling approaches and discuss their properties in the following section. 
For ease of presentation, we only consider resampling tests for $H_0^p:\{p=\frac{1}{2}\}$ in order to concentrate on one parameter of interest (i.e. on $p$) only. 
Nevertheless, the results directly carry over to construct resampling-based confidence intervals for $p$ and $w$, respectively.

\section{Resampling the Mann-Whitney statistic}
\label{sec:res}

Even in the continuous case, \cite{koziol2009concordance} pointed out that ``with small samples sizes a bootstrap approach might be preferable'' 
for approximating the unknown distribution of $V_n=V_n(p) = \sqrt{{n_1 n_2}/{n}} (\wh p - p)$. 
To this end, we consider different resampling methods, starting with Efron's classical bootstrap. 
Here the bootstrap sample is generated by drawing with replacement from the original data pairs; see \cite{efron81}. 
Large sample properties of the bootstrapped Kaplan-Meier process and extensions thereof have been analyzed
e.g. in \cite{akritas86}, \cite{lo86} and \cite{horvath87}.
Calculating the bootstrap version of $\wh p$ 
via bootstrapping for each sample group
and using their quantiles 
leads to a slightly improved control of the type-$I$ error probability 
in comparison to the asymptotic test~\eqref{eq:consistent_test}.
However, this way of bootstrapping results in a still too inaccurate behaviour
in terms of too large deviations from the $\alpha=5\%$ level (results not shown). 
This technique is typically improved by resampling procedures based on the pooled data
$\b Z = \{(Z_i, \eta_i): i=1,\dots,n\}$ given by
$$ (Z_i, \delta_i) = (X_{1i}, \delta_{1i}) \mathbf{1}\{ i \leq n_1 \} +  (X_{2 (i-n_1)}, \delta_{2 (i-n_1)}) \mathbf{1}\{ i > n_1 \}, \quad  i=1,\dots, n ;$$
see e.g. \cite{boos89}, \cite{janssen05monte} and \cite{neubert07} for empirical verifications for other functionals in this matter.
\cite{boos89} and \cite{konietschke14} also demonstrate that random permuting of and bootstrapping from pooled samples may yield to superior results, 
where the first has the additional advantage of 
leading to finitely exact testing procedures in case of $S_1=S_2$ and $G_1=G_2$. 
We investigate both techniques in more detail below.

\subsection{The pooled bootstrap}
\label{sec:ebs}

We independently draw  $n$ times  with replacement from the pooled data $\b Z$
to obtain the pooled bootstrap samples 
$\b Z_1^* = (Z_{1i}^*, \eta_{1i}^*)_{i=1}^{n_1}$ and $\b Z_2^* = (Z_{2i}^*, \eta_{2i}^*)_{i=1}^{n_2}$.
Denote the corresponding Kaplan-Meier estimators based on these bootstrap samples as $S_{1}^*$ and $S_{2}^*$.
These may also be regarded as the $n_j$ out of $n$ bootstrap versions of the Kaplan-Meier estimator $\wh S$ based on the pooled sample $\b Z$.
All in all, this results in the pooled bootstrap version
$p^* = - \int S_1^{*\pm} \d S_2^*$ of $\wh p$.
A suitable centering term for $p^*$ is based on the pooled Kaplan-Meier estimator
and is given by $ - \int \wh S^\pm \d \wh S = \frac12$.
Thus, we study the distribution of $V_n^* = \sqrt{\frac{n_1 n_2}{n}} (p^* - \frac12)$
for approximating the null distribution of $\sqrt{\frac{n_1 n_2}{n}} (\wh p - \frac12)$. 

To investigate the large sample behaviour of the bootstrap statistic $p^*$, first note that the pooled Kaplan-Meier estimator $\wh S$ is a functional 
of the empirical processes based on $\b X_1$ and $\b X_2$. 
Since this functional is Hadamard-differentiable with uniformly continuous linear derivative function, 
Donsker theorems for the empirical processes of $\b X_1$ and $\b X_2$ immediately carry over to the pooled Kaplan-Meier estimator. 
In case of continuously distributed event times this can also be seen by utilizing the usual martingale arguments of \cite{abgk93}, Section~IV.3.
In all convergence results stated below the c\`adl\`ag space $D[0,K]$ is always equipped with the $\sup$-norm; cf. \cite{vaart96}.
\begin{lemma}
\label{lem:pooled_kme_weak}
 Suppose that~\eqref{cond:liminfsup2} holds.
Then, as $n \rightarrow \infty$, we have $\sqrt{n} ( \wh S - S) \oDo U$ on $D[0,K]$, 
where $U$ is a zero-mean Gaussian process with covariance function 
\begin{align}
\label{eq:pooled_cov}
 \Gamma: (r,s) \mapsto S(r) S(s) \int_0^{r \wedge s}  \frac{\d \Lambda}{(1 - \Delta \Lambda)(\kappa S_{1-} G_{1-} + (1 - \kappa) S_{2-} G_{2-})}.
\end{align}
Here $ \d \Lambda = \frac{\kappa S_{1-} G_{1-}}{\kappa S_{1-} G_{1-} + (1 - \kappa) S_{2-} G_{2-}} \d \Lambda_1 + \frac{(1-\kappa) S_{2-} G_{2-} }{\kappa S_{1-} G_{1-} + (1 - \kappa) S_{2-} G_{2-}}\d \Lambda_2$ and $S(t) = \prodi_{(0,t]} (1- \d \Lambda)$.
\end{lemma}
A similar behaviour of the bootstrap counterpart is established in the following theorem.
Its proof relies on the $\delta$-method for the bootstrap from which the large sample properties of $V_n^*$ can be established as well.
\begin{thm}
 \label{thm:ebs kme}
 Suppose that~\eqref{cond:liminfsup2} holds.
    As $n \rightarrow \infty$ and given $\b Z$, we have conditional weak convergence on $D[0,K]$,
   $$ \sqrt{n_j} (S_{j}^* - \wh S) \oDo U, \quad j =1,2,$$ 
     in outer probability towards a Gaussian zero-mean process $U$ with covariance function $\Gamma$ given in~\eqref{eq:pooled_cov}.
\end{thm}
Since pooled sampling affects the covariance structure related to $V_{n}^*$ (see the Appendix for details), a studentization becomes mandatory. 
This is also in line with the general recommendation to bootstrap studentized statistics, see e.g. 
\cite{hall1991two}, \cite{janssen05monte} or \cite{delaigle2011robustness}. 
To this end, introduce the bootstrap variance estimator
\begin{align*}
  \sigma^{* 2} = \frac{n_2}{n} \int \int S_{1}^*(u) \Big[ n_1 \int_0^{u \wedge v} \frac{\d \Lambda_1^*}{(1 - \Delta \Lambda_1^*)Y_1^{*}} \Big] S_{1}^*(v)
  \d S_{2}^*(u) S_{2}^*(v)
  \\ + \frac{n_1}{n} \int \int S_{2}^*(u) \Big[ n_2 \int_0^{u \wedge v} \frac{\d \Lambda_2^*}{(1 - \Delta \Lambda_2^*)Y_2^{*}} \Big] S_{2}^*(v)
  \d S_{1}^*(u) S_{1}^*(v),
\end{align*}
where $N_j^*$ and $Y_j^*$ are the obvious bootstrap versions of the counting processes $N_j$ and $Y_j$
and $\d \Lambda_j^* = \d N_j^* / Y_j^*$ define the pooled bootstrap version of the Nelson-Aalen estimator, $j=1,2$.
We state our main result on the pooled bootstrap.
\begin{thm}
 \label{thm:ebs stud}
 Suppose that~\eqref{cond:liminfsup} holds.
 Then the studentized bootstrap statistic $T_{n}^*=V_{n}^*/ \sigma^{*}$ 
 always approximates the null distribution of $T_n(1/2)$ in outer probability, i.e. we have for any choice of $p$ and as $n_1 \wedge n_2 \rightarrow \infty$:
 \begin{equation}
  \sup_x \left| P_p(T_{n}^*\leq x|\vZ) - P_{1/2}\left(T_n(1/2)\leq x\right)\right| \oPo 0.
 \end{equation}
Moreover, denoting by $c_{n}^*(\alpha)$ the conditional $(1-\alpha)$-quantile of $T_{n}^*$ given $\b Z$,
it follows that $\varphi_{n}^* = \mathbf{1}\{T_n(1/2) > c_{n}^*(\alpha)\}$ 
is a consistent asymptotic level $\alpha$ test for $H_0^p:\{p=\frac{1}{2}\}$ 
against $H_1^p$ that is asymptotically equivalent to $\varphi_n$, 
i.e. we have $E(|\varphi_n^* - \varphi_n|)\to 0$.
\end{thm}

\subsection{Random permutation}
\label{sec:perm}

An alternative resampling technique to Efron's bootstrap is the permutation principle.
The idea is to randomly interchange the group association of all individuals while maintaining the original sample sizes.
The effect measure is then calculated anew based on the permuted samples.
A big advantage of permutation resampling over the pooled bootstrap 
is the finite exactness of inference procedures on the smaller null hypothesis
\begin{align}
\label{eq:H0_restr}
 H_0^{S,G}: \{ S_1 = S_2 \text{ and } G_1 = G_2 \} \subset H_0^p;
\end{align}
see e.g. \cite{neuhaus1993conditional} and \cite{brendel2014weighted} in case of testing $H_0^S$ and 
\cite{janssen97}, \cite{janssen1999},  \cite{neubert07}, \cite{ChungRomano2013, ChungRomano2013b},  \cite{pauly15} 
as well as \cite{pauly2016} in other situations without censoring.

Therefore, let $\pi: \Omega \rightarrow \mac{S}_n$ be independent of $\b Z $ 
and uniformly distributed on the symmetric group $\mac{S}_n$, the set of all permutations of $(1,\dots,n)$. 
The permuted samples are obtained as $\b Z_1^\pi = (Z_{\pi(i)}, \eta_{\pi(i)})_{i=1}^{n_1}$ and $\b Z_2^\pi = (Z_{\pi(i)}, \eta_{\pi(i)})_{i=n_1+1}^{n}$.
Plugging  the Kaplan-Meier estimators $S_1^{\pi}$ and $S_2^\pi$ based on these permuted samples
into the Wilcoxon functional leads to the permutation version
$$ {p}^\pi = - \int S_1^{\pi\pm} \d S_2^\pi $$
of $\wh p$.
The permutation sampling is equivalent to drawing without replacement from the pooled sample $\b Z$.
The following auxiliary result for the permutation version of the Kaplan-Meier estimator 
is analogous to Theorem~\ref{thm:ebs kme}, where now $D^2[0,K]$ is equipped with the $\max$-$\sup$-norm.
\begin{thm}
 \label{thm:perm kme}
 Suppose that~\eqref{cond:liminfsup2} holds. Then, as $n \rightarrow \infty$, the permutation versions of the Kaplan-Meier estimators conditionally converge on $D^2[0,K]$ in distribution
   $$ (\sqrt{n_1} (S_{1}^\pi - \wh S), \sqrt{n_2} (S_{2}^\pi - \wh S) )  \oDo (\sqrt{1 - \kappa} U, -\sqrt{\kappa} U) $$ 
    given $\b Z$ in outer probability. Here, $U$ is a zero-mean Gaussian process with covariance function $\Gamma$ given in \eqref{eq:pooled_cov}.
\end{thm}
It is shown in the Appendix, that $V_n^\pi = \sqrt{\frac{n_1 n_2}{n}} ( p^\pi - \frac12)$ does in general not possess the same limit distribution as $V_n$. 
Indeed, the limit variances may be different in general and again studentization of $V_n^\pi$ is necessary. This is achieved by utilizing the permutation version of 
$\widehat{\sigma}^2$:
\begin{align*}
  \sigma^{\pi 2} = \frac{n_2}{n} \int \int S_{1}^\pi(u) \Big[ n_1 \int_0^{u \wedge v} \frac{\d N_1^\pi}{(1-\frac{\Delta N_1^\pi}{Y_1^\pi}) Y_1^{\pi 2}} \Big] S_{1}^\pi(v)
  \d S_{2}^\pi(u) S_{2}^\pi(v)
  \\ + \frac{n_1}{n} \int \int S_{2}^\pi(u) \Big[ n_2 \int_0^{u \wedge v} \frac{\d N_2^\pi}{(1-\frac{\Delta N_2^\pi}{Y_2^\pi})Y_2^{\pi 2}} \Big] S_{2}^\pi(v)
  \d S_{1}^\pi(u) S_{1}^\pi(v)
\end{align*}
yielding the studentized permutation statistic $T_n^\pi=V_n^\pi/\sigma^{\pi}$. It is indeed the permutation version of $T_n(1/2)$ and also shares its asymptotic distribution as stated below.
\begin{thm}
 \label{thm:perm stud}
 Suppose that~\eqref{cond:liminfsup} holds.
 Then the studentized permutation statistic $T_n^\pi=V_n^\pi/\sigma^{\pi}$ always approximates the null distribution of $T_n(1/2)$ in outer probability, 
 i.e. we have for any choice of $p$  and as $n_1 \wedge n_2 \rightarrow \infty$:
 \begin{equation}
  \sup_x \left| P_p(T_n^\pi\leq x|\vZ) - P_{1/2}\left(T_n(1/2)\leq x\right)\right| \oPo 0.
 \end{equation}
Moreover, denoting by $c_{n}^\pi(\alpha)$ the conditional $(1-\alpha)$-quantile of $T_n^\pi$ given $\b Z$,  it follows that 
$
\varphi_{n}^\pi = \mathbf{1}\{T_n(1/2) > c_{n}^\pi(\alpha)\} 
$ possesses the same asymptotic properties as $\varphi_n^*$ in Theorem~\ref{thm:ebs stud}.
Furthermore, $\varphi_n^\pi$ is even a finitely exact level $\alpha$ test under $H_0^{S,G}$.
\end{thm}

\section{Finite sample properties}
\label{sec:simus}

In this section we study the finite-sample properties of the proposed approximations. In particular, we compare the actual coverage probability 
of the asymptotic two-sided confidence interval $I_n$ given in \eqref{eq:consistent_cis} with that of the corresponding bootstrap and permutation confidence intervals 
$$
I_n^* = \left[\wh p \ \mp \ \frac{c_n^*(\alpha/2)\wh \sigma \sqrt{n}}{\sqrt{n_1n_2}}  \right] \quad\text{ and } \quad I_n^\pi = \left[\wh p \ \mp \ \frac{c_n^\pi(\alpha/2)\wh \sigma \sqrt{n}}{\sqrt{n_1n_2}}  \right],
$$
respectively.
To this end, the following distribution functions $F_1$ and $F_2$, frequently occuring in the survival context, have been chosen in our simulation study:
\begin{itemize}
 \item[(1)] Group 1: Exponential distribution with mean $1/2$, i.e. $F_1=Exp(1/2)$. \\
	    Group 2: Mixture of two exponential distributions: \\
	    $F_2= \frac{1}{3}Exp(1/1.27) + \frac{2}{3}Exp(1/2.5)$.\\ 
	    End-of-study time: $K \approx 1.6024$ such that $p \approx 1/2$.
 \item[(2)] Group 1: Weibull distribution with scale parameter 1.65 and shape parameter 0.9. \\
	    Group 2: Standard lognormal distribution. \\
	    End-of-study time: $K \approx 1.7646$ such that $p \approx 1/2$.
 \item[(3)] Groups 1 and 2: Equal Weibull distributions with scale parameters 1 and shape parameters 1.5. \\
	    End-of-study time: $K = 2$ such that $p = 1/2$.
\end{itemize}
Censoring is realized using i.i.d. exponentially distributed censoring variables $C_{ji}$ with parameters chosen 
such that the (simulated) censoring probability (after truncation at $K$) for each of both sample groups belongs to the following ranges:
\begin{itemize}
 \item Strong censoring: Censoring percentages between $40.97$ and $43.6$ per cent,
 \item Moderate censoring: Censoring percentages between $21.19$ and $26.39$ per cent,
 \item No censoring: Censoring percentages set to zero.
\end{itemize}
The sample sizes range over $n_1 = n_2 \in \{10,15,20,25,30\}$ as well as $n_2 = 2 n_1 \in \{20,30, 40,50, 60\}$.
Simulating  10,000 individuals each,
the approximate proportions of observations greater than $K$ are given in Table~\ref{tab:propK}.
%
%
%
\begin{table}[H]
\centering
\begin{tabular}{|c|c|c|c|} \hline
\mbox{ } & strong & moderate & none \\
set-up / censoring & Group I/II & Group I/II & Group I/II \\
\hline
(1) & 0.36 / 0.32 	& 1.19 / 1.35 	& 4.12 / 4.26 	\\
(2) & 12.16 / 10.44 	& 22.3 / 17.43 	& 34.26 / 28.37 \\
(3) & 1.54 	& 3.28 		& 6.02		\\ \hline
\end{tabular}
\caption{Simulated percentages of observations greater than the respective $K$.
In set-ups (1) and (2) the proportions for both sample groups are separated by a `` \ / \ ''. }
\label{tab:propK}
\end{table}
The pre-specified nominal level is $1-\alpha = 95\%$.
Each simulation was carried out using $N=10,000$ independent tests, each with $B=1,999$ resampling steps
in \verb=R= version 3.2.3 (\citealt{Rteam}).
All Kaplan-Meier estimators were calculated using the \verb=R= package {\it etm} by \cite{etm}.
In comparison to the pooled bootstrap confidence interval $I_n^*$, 
the permutation-based confidence interval $I_n^\pi$ provides even finitely exact coverage probabilities if the restricted null hypothesis $H_0^{S,G}$ given in~\eqref{eq:H0_restr} is true.

The simulation results for all scenarios are summarized in Tables~\ref{tab:gl} and \ref{tab:ug}. 
Starting with the balanced setting ($n_1=n_2$) it can be readily seen that in case of no censoring,
both, the asymptotic and permutation confidence intervals $I_n$ and $I_n^\pi$ show very accurate coverages; even for small samples. 

\begin{table}[H]
\centering
\begin{tabular}{|c|c|ccc|ccc|ccc|}
  \hline
 & censoring & \multicolumn{3}{c}{strong}  & \multicolumn{3}{c}{moderate} & \multicolumn{3}{c|}{none}  \\ \hline
  set-up &  $n_1$ & $I_n$ & $I_n^*$ & $I_n^\pi$ & $I_n$ & $I_n^*$ & $I_n^\pi$ & $I_n$ & $I_n^*$ & $I_n^\pi$ \\ \hline 
  (1) & 10, 10 & 90.63 & 95.70 & 94.95 & 93.49 & 94.63 & 95.08 & 95.02 & 93.80 & 94.63 \\ 
   & 15, 15 & 92.69 & 95.48 & 94.99 & 94.26 & 94.59 & 95.06 & 94.76 & 93.91 & 94.46 \\ 
   & 20, 20 & 93.08 & 95.83 & 94.65 & 94.00 & 94.20 & 94.63 & 95.56 & 93.84 & 94.98 \\ 
   & 25, 25 & 93.64 & 95.37 & 94.38 & 94.41 & 94.18 & 94.70 & 94.95 & 93.45 & 94.66 \\ 
   & 30, 30 & 93.96 & 94.91 & 94.31 & 94.73 & 94.12 & 94.27 & 95.34 & 93.80 & 94.37 \\  \hline
  (2) & 10, 10 & 91.73 & 95.32 & 94.85 & 93.97 & 94.26 & 95.26 & 95.48 & 93.68 & 94.97 \\ 
   & 15, 15 & 93.26 & 94.73 & 95.08 & 94.76 & 94.66 & 94.69 & 95.48 & 94.11 & 95.18 \\ 
   & 20, 20 & 93.70 & 94.80 & 95.00 & 94.80 & 94.47 & 95.05 & 95.48 & 94.17 & 95.21 \\ 
   & 25, 25 & 94.06 & 94.96 & 94.95 & 94.87 & 94.41 & 95.13 & 95.66 & 94.09 & 94.78 \\ 
   & 30, 30 & 94.28 & 94.52 & 94.69 & 95.22 & 94.78 & 94.86 & 95.64 & 94.44 & 94.74 \\  \hline
  (3) & 10, 10 & 89.88 & 95.33 & 94.64 & 92.86 & 94.65 & 94.92 & 94.86 & 94.30 & 94.97 \\ 
   & 15, 15 & 92.05 & 95.60 & 94.82 & 93.85 & 94.58 & 95.09 & 95.57 & 94.06 & 95.28 \\ 
   & 20, 20 & 92.98 & 95.05 & 95.32 & 94.33 & 94.26 & 95.03 & 95.22 & 94.36 & 95.25 \\ 
   & 25, 25 & 93.32 & 95.12 & 94.78 & 94.07 & 94.32 & 95.11 & 95.34 & 94.38 & 95.04 \\ 
   & 30, 30 & 93.90 & 95.09 & 94.96 & 94.63 & 94.60 & 95.09 & 94.95 & 94.05 & 94.64 \\ 
   \hline
\end{tabular}
\caption{Simulated coverage probabilities (in \%) of two-sided asymptotic $95\%$ confidence intervals for $p = 0.5$ and equal sample sizes $n_1=n_2$.} 
\label{tab:gl}
\end{table}

This is in line with previous findings of $\alpha$-level control of rank-based tests for $H_0^p$ (e.g. \citealt{neubert07} or \citealt{pauly2016}). 
In comparison the bootstrap based confidence intervals 
$I_n^*$ are slightly liberal; especially in the exponential set-up (1) with coverages between $93.45\%$ and $93.91\%$.  
 In case of moderate and strong censoring, however, the behaviour of the asymptotic and bootstrap confidence intervals changes: 
The asymptotic procedure $I_n$ based on normal quantiles gets liberal while the bootstrap procedure is more or less accurate. 
The liberality of $I_n$ is particularly apparent for the smallest total sample sizes of $n=20$ (with coverages between $89.88\%$ and $91.73\%$) but 
also remains present with increasing sample sizes.
In comparison, the permutation interval $I_n^\pi$ achieved very accurate results in all set-ups. 
Even in the presence of strong censoring and with unequal distributions (set-ups (1) and (2)) 
the simulated coverage probabilities are always very close to the nominal level of 95\%.

\begin{table}[H]
\centering
\begin{tabular}{|c|c|ccc|ccc|ccc|}
  \hline
 & censoring & \multicolumn{3}{c}{strong}  & \multicolumn{3}{c}{moderate} & \multicolumn{3}{c|}{none}  \\ \hline
  set-up &  $n_1$ & $I_n$ & $I_n^*$ & $I_n^\pi$ & $I_n$ & $I_n^*$ & $I_n^\pi$ & $I_n$ & $I_n^*$ & $I_n^\pi$ \\ \hline 
  (1) & 10, 20 & 91.66 & 95.21 & 95.13 & 93.08 & 94.09 & 94.91 & 94.20 & 93.71 & 94.48 \\ 
   & 15, 30 & 92.41 & 95.35 & 94.59 & 93.77 & 94.16 & 95.00 & 94.34 & 93.86 & 94.30 \\ 
   & 20, 40 & 93.62 & 94.77 & 94.53 & 93.93 & 93.93 & 94.63 & 94.50 & 93.60 & 94.51 \\ 
   & 25, 50 & 93.07 & 94.55 & 94.37 & 93.98 & 94.06 & 94.01 & 94.46 & 93.96 & 94.43 \\ 
   & 30, 60 & 92.94 & 94.42 & 94.48 & 93.97 & 94.28 & 94.46 & 94.83 & 93.27 & 94.35 \\  \hline
  (2) & 10, 20 & 92.02 & 95.21 & 95.11 & 93.43 & 94.49 & 95.09 & 94.93 & 94.33 & 95.03 \\ 
   & 15, 30 & 92.94 & 95.30 & 95.37 & 94.57 & 94.90 & 95.16 & 95.29 & 94.60 & 95.32 \\ 
   & 20, 40 & 94.01 & 95.45 & 95.29 & 94.45 & 94.75 & 95.19 & 95.46 & 95.67 & 95.98 \\ 
   & 25, 50 & 93.94 & 95.21 & 94.93 & 95.14 & 94.84 & 95.11 & 95.72 & 95.15 & 95.27 \\ 
   & 30, 60 & 94.24 & 94.97 & 94.93 & 94.73 & 95.30 & 95.00 & 95.74 & 94.81 & 96.05 \\  \hline
  (3) & 10, 20 & 90.99 & 95.46 & 95.14 & 92.99 & 94.79 & 94.88 & 94.61 & 94.27 & 94.63 \\ 
   & 15, 30 & 92.18 & 95.07 & 95.23 & 93.58 & 94.90 & 95.12 & 94.61 & 94.00 & 94.81 \\ 
   & 20, 40 & 93.35 & 95.27 & 95.04 & 94.02 & 94.80 & 94.72 & 94.81 & 94.82 & 94.72 \\ 
   & 25, 50 & 93.18 & 95.51 & 95.14 & 94.78 & 94.80 & 94.58 & 94.84 & 94.59 & 94.97 \\ 
   & 30, 60 & 94.03 & 94.82 & 95.01 & 94.45 & 95.02 & 94.82 & 94.90 & 94.24 & 95.04 \\ 
   \hline
\end{tabular}
\caption{Simulated coverage probabilities (in \%) of two-sided asymptotic $95\%$ confidence intervals for $p = 0.5$ and unequal sample sizes $2n_1=n_2$.} 
\label{tab:ug}
\end{table}

From the results for the unbalanced case $2n_1=n_2$, shown in Table~\ref{tab:ug}, we can draw the same conclusion. 
In particular, the permutation approach seems to be the most promising since it does not only show the best coverages probabilities 
but has the additional advantage of being finitely exact in case of 
equal censoring and survival distributions, i.e. for $F_1=F_2$ (implying $p=1/2$) and $G_1=G_2$. 

All in all, the permutation procedure can be generally recommended. In the censored case the bootstrap procedure shows a similar coverage (with a minor liberality for moderate censoring) 
but does not possess the nice exactness property under $H_0^{S,G}$. The asymptotic procedure can only be recommended for larger sample sizes if no or only slight censoring is apparent.

%
%
%
%

%
%
%
%
%
\section{Application to a data example}
\label{sec:example}
%
%
%
%

To illustrate the practical applicability of our novel approaches we reconsider a data-set containing survival times of tongue cancer patients, 
cf.  \cite{klein03}. 
This data-set is freely available in the \verb=R= package {\it KMsurv} via \verb=data(tongue)=.
It contains 80 patients of which $n_1=52$ are suffering from an aneuploid tongue cancer tumor (group 1) and $n_2=28$ are suffering from a diploid tumor (group 2).
Observation of 21 patients in group 1 and of six patients in group 2 have been right-censored, for all others the time of death has been recorded.
Thus, the corresponding censoring proportions are intermediate between the ``strong'' and ``moderate'' scenarios of Section~\ref{sec:simus}. 
Note that the data set actually contains ties: among the uncensored survival times, there are 27 different times of death in the first group, 20 different in the second group, and 39 different in the pooled sample.
There are three individuals in group 1 with censoring time exceeding the greatest recorded time of death in this group; 
for group 2 there is one such individual.
As a reasonable value for restricting the time interval,
we may thus choose $K=200$ weeks which still precedes all just mentioned censoring times.

\begin{figure}[ht]
\centering
 \includegraphics[width=1.03\textwidth, height=0.55\textwidth]{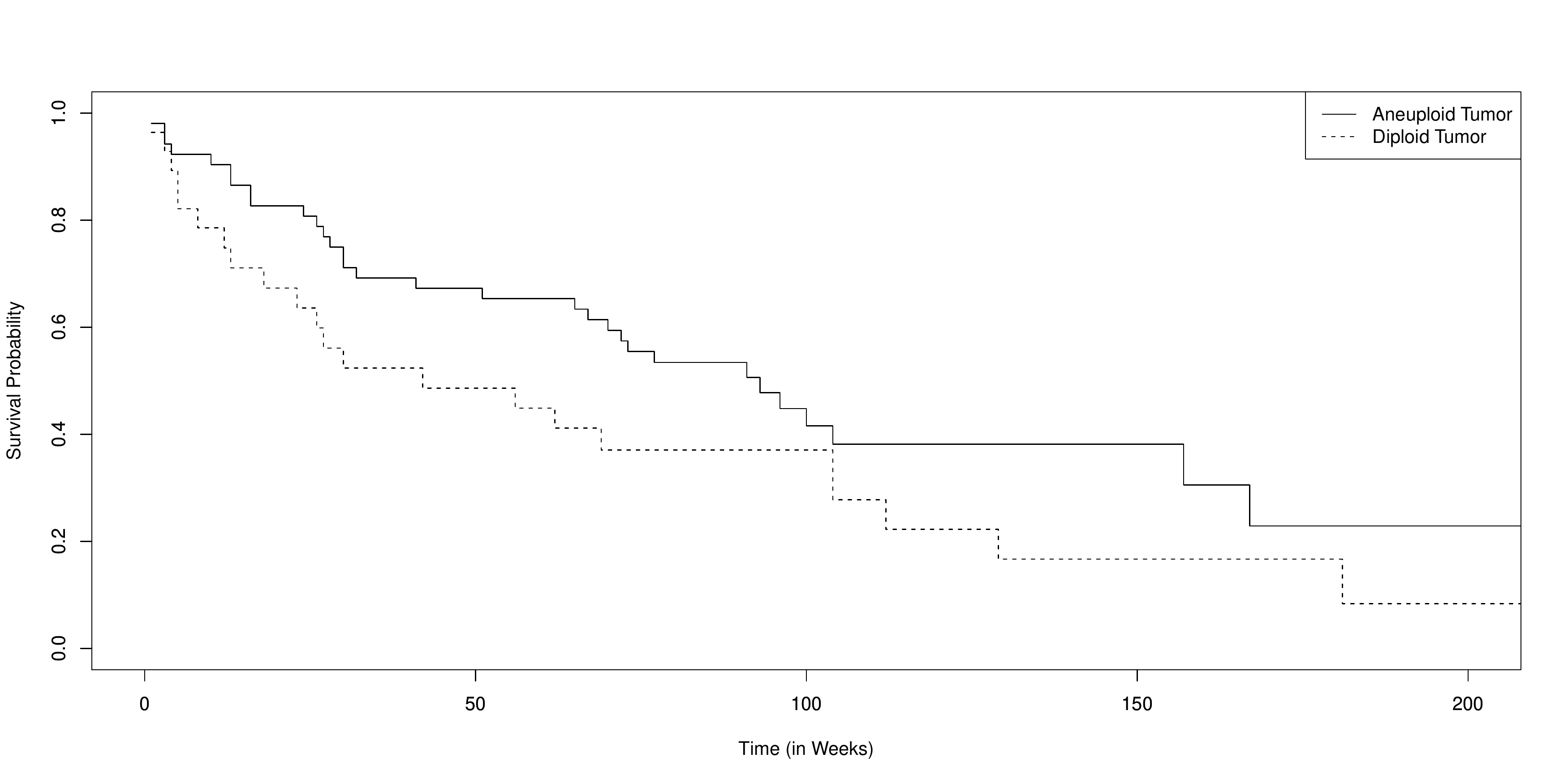}
 \caption{Kaplan-Meier estimators for patients with diploid tumor (- - -) and aneuploid tumor (-----).}
 \label{fig:kmes}
\end{figure}

The Kaplan-Meier estimators correspoding to both recorded groups are plotted in Figure~\ref{fig:kmes}.
It shows that the aneuploid Kaplan-Meier curve is always above the Kaplan-Meier curve of the diploid group.
We would therefore like to examine whether this gap already yields significant results concerning the probability of concordance.
%
%
%

The data evaluation resulted in a point estimate $\wh p \approx 0.6148$ indicating 
a slightly larger survival probability of the aneuploid group in comparison to the diploid. 
To infer this,  
the one- and two-sided $95\%$ confidence intervals (based on normal, bootstrap and permutation quantiles) for the probability $p$ 
that a randomly chosen individual with an aneuploid tumor survives longer than a corresponding patient with diploid tumor are given in Table~\ref{tab:CIs}.
These have been calculated using the asymptotic normal quantile as well as $B=9,999$ resampling iterations for each of the bootstrap and the permutation technique.

\begin{table}[H]
\centering
\begin{tabular}{|c|c|c|} \hline
\mbox{} & \multicolumn{2}{c|}{confidence interval}\\
method & two-sided & one-sided\\\hline
asymptotic & [0.475, 0.755] & [0.497, 1.000]\\
bootstrap& [0.457, 0.772] & [0.507, 1.000]\\
permutation& [0.464, 0.766] & [0.506, 1.000]
\\ \hline
\end{tabular}
\caption{One- and two-sided $95\%$ confidence intervals based on normal, bootstrap and permutation quantiles for the Mann-Whitney effect in the tongue cancer data. 
In all cases the point estimate is $\wh p \approx 0.6148$.}
\label{tab:CIs}
\end{table}

By inverting these confidence intervals it can be readily seen that the two-sided null hypothesis $H_0^p:\{p=1/2\}$ cannot be rejected by any procedure 
since $p=1/2$ is contained in all two-sided intervals. However, if we were only interested in detecting an effect in favor of the 
aneuploid group, we have to consider the corresponding one-sided tests to 
avoid possible directional errors. In particular, the results for testing the one-sided hypothesis $H_{0,\leq}^p:\{p\leq 1/2\}$ are borderline: It can be rejected by the resampling approaches 
at level $5\%$ but lies close to the confidence limit of the asymptotic interval. 
Here, a slightly larger data-set might have caused a different decision. We note that multiplicity issues have not been taken into account.

\section{Summary and discussion}
\label{sec:dis}

{
In this article, novel inference procedures for the Mann-Whitney effect $p$ and the win ratio $w$ are introduced {\color{black}both of} which are meaningful and well-established effect measures (especially
in biometry and survival analysis). 
In comparison to the usual survival hypothesis $H_0^S:\{S_1=S_2\}$, 
we were the first who particularly developed asymptotic confidence intervals for $p$ and $w$ as well as tests for the more interesting 
composite null hypothesis $H_0^p:\{ p= 1/2\}$ in the two-sample survival model with right-censored data. 
By utilizing normalized Kaplan-Meier estimates these can even be constructed for discontinuous{\color{black}ly distributed} survival times that may be subject to 
independent right-censoring. 
Applying empirical process theory we showed that point estimates of $p$ and $w$ are asymptotically normal.
By introducing novel variance estimates, this leads to asymptotic inference procedures based on normal quantiles. 
To improve their finite sample performance, bootstrap and permutation approaches have been considered and shown to maintain the same asymptotic properties.
In our simulation study it could be seen that the  proposed permutation procedure considerably improves the finite sample performance of our procedure. 
Moreover, it is even finitely exact if data is exchangeable (i.e., whenever both survival and both censoring distributions are equal) 
and can thus be recommended as the method of choice. 
In the special continuous situation with complete observations a similar result has been recently proven in \cite{chung2016asymptotically}.



Note, that the proposed method can also be applied in the `winner-loser' set-ups 
considered in \cite{pocock2012win} or \cite{wang2016win}, where 
now even the neglected ties can be taken into account. We plan to do this in the near future. 
Moreover, extensions of the proposed techniques to other models such as multiple samples  and multivariate \nocite{davidov2013linear} or specific paired designs (e.g. measurements before and after treatment) 
will also be considered in a forthcoming paper.
}
 
\appendix 

\vspace{1cm}

\noindent {\Large {\bf Appendix}}

\vspace{-1cm}

\section{Proofs}
\label{sec:proofs}

\noindent{\it Proof of Theorem~\ref{thm:weak}}:
Integration by parts shows that
\begin{align*}
p & = - \int S_1^{\pm} \d S_2 
= - \frac12 \int S_1 \d S_2 + \frac12 \int S_2 \d S_1 - \frac12 S_1(K) S_2(K) + \frac12 S_1(0) S_2(0) \\
& = - \frac12 \int S_1 \d S_2 + \frac12 \int S_2 \d S_1 + \frac12 \\
& = - \frac12 \int_{[0,K)} S_1 \d S_2 + \frac12 \int_{[0,K)} S_2 \d S_1 + \frac12 =: \phi(S_1,S_2).
\end{align*}
Here it is important to exclude the right boundary $K$ of the support of $S_1$ and $S_2$
since otherwise theorems on the weak convergence of the Kaplan-Meier estimator on the whole support of the survival function would be required;
see \cite{gill83} and \cite{ying89} for such statements as well as \cite{dobler16_BS_KME_WL} for a bootstrap version.
Each integral is a Hadamard-differentiable functional of $(S_1,S_2)$ tangentially to 
\begin{align*}
 \mac H^2 := \{ (h_1, h_2) \in D^2[0,K] : \ &  h_1^\pm \in \mac L_1|_{[0,K]}(S_2), h_2^\pm \in \mac L_1|_{[0,K]}(S_1), \\
   & h_1(0) = h_2(0),  h_1(K) = h_2(K) = 0 \};
\end{align*}
see e.g. \cite{vaart96}, Lemma~3.9.17.
Note that $S_1, S_2$ are monotone and thus of bounded variation.
The continuous, linear Hadamard-derivative is given by
$$ \d \phi_{(S_1,S_2)} \cdot (h_1, h_2) = \frac12 \Big[ - \int_{[0,K)} S_1 \d h_2 - \int_{[0,K)} h_1 \d S_2 + \int_{[0,K)} S_2 \d h_1 + \int_{[0,K)} h_2 \d S_1 \Big], $$
integrals with respect to $h_1,h_2$ defined via integration by parts.
Because the integrands and integrators have no mass in $K$, this derivative can be further simplified to
\begin{align}
\label{eq:hadamard_der_wilcox}
  \d \phi_{(S_1,S_2)} \cdot (h_1,h_2) = - \int h_1^\pm \d S_2 + \int h_2^\pm \d S_1,
\end{align}
again via integration by parts.
By the functional $\delta$-method (Section~3.9 in \citealt{vaart96}), we thus conclude that
\begin{align*}
 & V_n(p) =  \sqrt{\frac{n_1 n_2}{n}} [ \phi(\wh S_1, \wh S_2) - \phi(S_1, S_2)]  \\
 & = \d \phi_{(S_1,S_2)} \cdot \sqrt{\frac{n_1 n_2}{n}} (\wh S_1 - S_1, \wh S_2 - S_2) + o_p(1) \\
 & = \sqrt{\frac{n_1}{n}} \int \sqrt{n_2} (\wh S_2 - S_2)^\pm \d S_1
   - \sqrt{\frac{n_2}{n}} \int \sqrt{n_1} (\wh S_1 - S_1)^\pm \d S_2 + o_p(1)\\
 & \oDo \sqrt{\kappa} \int U_2^{\pm} \d S_1
  + \sqrt{1-\kappa} \int U_1^{\pm} \d S_2 =: Z
\end{align*}
as $n \rightarrow \infty$, where $U_1$ and $U_2$ are independent, $D[0,K]$-valued, zero-mean Gaussian processes with covariance functions
$$ \Gamma_j : (r,s) \mapsto 1\{r,s < K \} S_j(r) S_j(s) \int_0^{r \wedge s} \frac{\d \Lambda_j}{(1-\Delta \Lambda_j) S_{j-} G_{j-}}, \quad j=1,2, \ 0 \leq r,s \leq K $$
and with $U_1(K) = U_2(K) \equiv 0$;
see Example~3.9.31 in \cite{vaart96}.
Thus, $E(Z)=0$ and
\begin{align*} \sigma^2 = var(Z) = \kappa \int \int \Gamma_2^{\pm \pm} \d S_1 \d S_1 
+ (1-\kappa) \int\int \Gamma_1^{\pm \pm} \d S_2 \d S_2 
\end{align*}
where again $\Gamma_j^{\pm \pm}(r,s) = \frac14 [ \Gamma_j(r,s) + \Gamma_j(r,s-) + \Gamma_j(r-,s) + \Gamma_j(r-,s-) ], \ j=1,2$.
$ \hfill \Box$

\noindent{\it Proof of Lemma~\ref{lem:cons_var_est}}:
The consistency $\wh \sigma^2 \oPo \sigma^2$ follows from the consistency of the Kaplan-Meier estimator,
the continuity of $D^2[0,K]\ni (f,g)  \mapsto \int f \d g$ in all functions $g$ of bounded variation (see the proof of Lemma~3 in \citealt{gill89}),
as well as from the continuity of $D[0,K]\ni f \mapsto \frac1f$ in all functions that are bounded away from zero.
$ \hfill \Box$

\noindent{\it Proof of Lemma~\ref{lem:pooled_kme_weak}}:
  This statement is proven basically in the same way as the large-sample properties of the sample-specific Kaplan-Meier estimators, i.e., applying the functional $\delta$-method
  (Theorem~3.9.5 in \citealt{vaart96})
  to the (pooled) empirical process of $\b Z$, indexed by $\mac F = \{ 1\{ \ \cdot_1 \ \leq z , \ \cdot_2 = 1 \}, \ 1\{ \ \cdot_1 \ \geq z \}: \ z \in [0,K) \}$.
  The only difference is the limit of the pooled Nelson-Aalen estimator:
  Writing $N = N_1 + N_2$ and $Y = Y_1 + Y_2$ and letting $t \in [0,K]$,
  we have
  \begin{align*}
    \wh \Lambda(t) = & \int_0^t \frac{\d N}{Y} = 
  \frac{n_1}{n} \int_0^t \frac{Y_1}{n_1} \frac{n}{Y} \frac{\d N_1}{Y_1} 
  + \frac{n_2}{n} \int_0^t \frac{Y_2}{n_2} \frac{n}{Y} \frac{\d N_2}{Y_2} \\
   \oPo  \kappa & \int_0^t  \frac{S_{1-} G_{1-}}{\kappa S_{1-} G_{1-} + (1-\kappa) S_{2-} G_{2-}} \d \Lambda_1
  + (1-\kappa) \int_0^t \frac{S_{2-} G_{2-}}{\kappa S_{1-} G_{1-} + (1-\kappa) S_{2-} G_{2-}} \d \Lambda_2 
  .
  \end{align*}
  Thus, substituting this quantity for the cumulative hazard function in the asymptotic covariance function of the Kaplan-Meier estimator, the proof is complete.
$ \hfill \Box$

\noindent{\it Proof of Theorem~\ref{thm:ebs kme}}:
  This is simply an application of the $\delta$-method for the bootstrap (Theorem~3.9.11 in \citealt{vaart96}) to the result of Lemma~\ref{lem:pooled_kme_weak}
  and the two-sample bootstrap Donsker Theorem~3.7.6 in \cite{vaart96}.
  Note, that we employ a Donsker theorem for 
  $$\sqrt{\frac{n_1 n_2}{n}} \Big(\mb P_n^\b Z - \frac{n_1}{n} P^{(X_1,\delta_1)} - \frac{n_2}{n} P^{(X_n,\delta_n)}\Big),$$
  $\mb P_n^\b Z$ being the empirical process of the pooled sample $\b Z$.
  This yields the same limit distribution of $\sqrt{n_j}(S_j^* - \wh S)$ as in Lemma~\ref{lem:pooled_kme_weak} by the uniform Hadamard-differentiability of all involved functionals, 
  concluding the proof.$ \hfill \Box$

\noindent{\it Proof of Theorem~\ref{thm:ebs stud}}:
  This is just another application of the functional $\delta$-method for the bootstrap, applied to the Wilcoxon statistic $\phi$ and the intermediate Theorem~\ref{thm:ebs kme}.
  To see the consistency of the bootstrap variance estimator,
  combine the Glivenko-Cantelli theorems for the bootstrap empirical processes of $\b Z_1^*$ and $\b Z_2^*$ with the continuous mapping theorem.
  Finally, apply Slutsky's lemma; e.g. Example~1.4.7 in \cite{vaart96}.
$ \hfill \Box$

\noindent{\it Proof of Theorem~\ref{thm:perm kme}}:
  This follows in the same way as Theorem~\ref{thm:ebs kme},
  but with the asymptotic covariance function determined by the Donsker theorem for the permutation empirical process (cf. Theorem~3.7.1 in \citealt{vaart96}) 
  and the observation that 
  $$\mb P_{n_1}^{\b X_1} - \mb P_n^\b Z \ = \ \frac{n_2}{n} (\mb P_{n_1}^{\b X_1} - \mb P_{n_2}^{\b X_2}) \ = \ - \frac{n_2}{n_1} (\mb P_{n_2}^{\b X_2} - \mb P_n^\b Z ).$$
  Here $\mb P_{n_1}^{\b X_1}$ and $\mb P_{n_2}^{\b X_2}$ are the empirical processes based on $\b X_1$ and $\b X_2$, respectively; cf. p. 361 in \cite{vaart96}.
$ \hfill \Box$

\noindent{\it Proof of Theorem~\ref{thm:perm stud}}:
  This theorem is proven in the same way as Theorem~\ref{thm:ebs stud}, using the intermediate result of Theorem~\ref{thm:perm kme}.
  However, the asymptotic variance of $V_n^\pi$ requires some attention:
  By the functional $\delta$-method $V_n^\pi$ is asymptotically equivalent to the image of 
  $$\sqrt{\frac{n_1 n_2}{n}} (\mb P_{n_1}^{\b Z_1^\pi} - \mb P_{n}^{\b Z}, \mb P_{n_2}^{\b Z_2^\pi} - \mb P_{n}^{\b Z}) \oDo ((1-\kappa) U, -\kappa U)$$ 
  under a continuous linear map,
  where $\mb P_{n_j}^{\b Z_j^\pi}$ is the $j$th permutation empirical process, $j=1,2,$
  and the Gaussian process $U$ is given by Theorem~\ref{thm:perm kme}.
  But this map (i.e., the Hadamard derivative) also subtracts both components of the previous display, cf.~\eqref{eq:hadamard_der_wilcox}.
  Hence, $T_n^\pi$ has the same large sample behaviour as its bootstrap counterpart $T_n^*$ of Theorem~\ref{thm:ebs stud}.
  The finite exactness under exchangeable data is a well-known property of permutation tests, see e.g. \cite{janssen2007most}.
$ \hfill \Box$

%
%
%
%
%







\renewcommand{\refname}{References}
\bibliography{literatur}
\bibliographystyle{plainnat}

\newpage

\end{document}